\input amstex
\documentstyle{amsppt}
\input bull-ppt

\ifx\undefined\rom\define\rom#1{{\rm #1}}\fi

\define\CPtwo{\Bbb C\Bbb P^2}
\define\CPtwobar{\rlap{$\bar{\phantom{\Bbb C\Bbb 
P}}$}\CPtwo}

\define\SO{\mathop{SO}\nolimits}

\topmatter
\cvol{30}
\cvolyear{1994}
\cmonth{April}
\cyear{1994}
\cvolno{2}
\cpgs{215-221}
\title
Recurrence relations and asymptotics\\ for
four-manifold invariants%
\endtitle
\ratitle
\shorttitle{Recurrence Relations for Four-Manifold 
Invariants}
\author	P. B. Kronheimer and T. S. Mrowka\endauthor
\address Mathematical Institute, 24--29 St.~Giles, Oxford 
OX1 3LB, England\endaddress
\ml kronheim\@maths.oxford.ac.uk\endml
\address California Institute of Technology, Pasadena, 
California 94720\endaddress
\ml mrowka\@cco.caltech.edu\endml
\subjclass Primary 57R57, 57N13, 57R55; Secondary 57R40, 14J60\endsubjclass
\date May 19, 1993\enddate
\abstract The polynomial invariants $q_d$ for a large  
class of smooth
4-manifolds are shown to satisfy universal 
relations.  The
relations reflect the possible genera of embedded surfaces
in the 4-manifold and lead to a structure theorem for the
polynomials.  As an application, one can read off a lower
bound for the genera of embedded surfaces from the
asymptotics of $q_d$ for large $d$.  The relations are
proved using moduli spaces of singular instantons.
\endabstract
\endtopmatter

\document

\heading1. Introduction\endheading

Donaldson's polynomial invariants \cite{3} are invariants 
of smooth, oriented,
compact 4-manifolds $X$ without boundary.  They and their 
close cousins are
the only useful invariants we have at present to 
distinguish different
4-manifolds of the same homotopy type.  We recall that 
when $X$ is simply
connected and $b^+(X)$ is odd and not less than three, the 
polynomial
invariants take the form of homogeneous polynomial functions
 $$
q_d : H_2(X, \Bbb R) \to \Bbb R.
 $$
Here $b^+$ is the dimension of a maximal positive subspace 
for the
intersection form on $H_2(X)$.  The index $d$ is the 
degree of $q_d$, and the
invariants are defined for all non-negative degrees
$d\equiv{1\over2}(b^++1)$ mod $4$. Formally, one can 
regard $q_d$ as being
given by the formula  
 $$
q_d(h) = \bigl\langle\, \mu(h)^d,\ [M_d]\,\bigr\rangle,
 $$
where $M_d$ is the instanton moduli space of dimension 
$2d$ (depending on a
choice of a Riemannian metric on $X$) and $\mu$ is a 
natural map from
$H_2(X)$ to $H^2(M_d)$. Because $M_d$ is non-compact in 
general, this
pairing needs to be correctly interpreted before it can be 
regarded as
well defined: a recipe for evaluation was given in 
\cite{3}, subject to the
constraint $2d>3(b^++1)$; and by various devices the 
construction has
since been extended \cite{6, 11} so as to remove this 
restriction.

Although the polynomial invariants have been a powerful 
tool in enlarging
our understanding of smooth 4-manifolds, their nature has 
remained
mysterious. The first applications rested on a small 
number of rather
general properties, such as the vanishing theorem for 
connected sums proved
in \cite{3}.  Other results have arisen from calculations 
of particular
values of $q_d$ for carefully chosen 4-manifolds; the 
elliptic surfaces are
a good example \cite{1, 5, 12, 13}.  Until now, however, 
there has remained only
one 4-manifold, the $K3$ surface, for which the invariants 
$q_d$ are
entirely known (apart, that is, from cases to which the 
vanishing theorem
applies).  The isolated values computed in other examples 
have remained a
short but slowly growing list.

The results which we announce here are not in the line of 
further particular
calculations.  Rather, we prove that there is a family of 
universal linear
relations which constrain the values of $q_d$, for a large 
class of
4-manifolds (the manifolds of {\it simple type\/}; see 
section~2).  It turns
out that the structure revealed by these relations is best 
expressed by
combining all the polynomials $q_d$ into one analytic 
function on
$H_2(X,\Bbb R)$: our main result (Theorem~1) then 
describes the structure of
this analytic function in terms of a distinguished 
collection of
2-dimensional cohomology classes, which we call the {\it 
basic classes\/} on
$X$.  Theorem~2 shows that the basic classes determine a 
lower bound on the
genus of embedded surfaces in the 4-manifold; the bound is 
reminiscent of
the way in which the canonical class of a complex surface 
determines the
genus of embedded complex curves.  

Section~4 explains the origin of the universal relations 
and gives some
pointers to the proofs of the theorems.  Our main tools 
are the moduli spaces
of instantons with singularities in codimension 2, 
introduced in \cite{9, 10}.

\heading2. Statement of results\endheading

Before stating our results, we must introduce the simple 
type condition (see
\cite{8}).  In addition to the 2-dimensional cohomology 
classes $\mu(h)$,
the moduli spaces $M_d$ also carry a 4-dimensional class 
$\nu$, and one can
define a larger class of invariants of $X$ by evaluating 
products such as
$\mu(h)^a \nu^b$ on the fundamental class $[M_d]$, 
whenever $a+2b=d$; see
\cite{4} for example. We adopt
the convention that $\nu$ is $-{1\over4}p_1(\Bbb E)$, 
where $\Bbb E$ is the
principal $\SO(3)$ bundle associated to the base-point 
fibration.
We say that $X$ has simple type if the 4-dimensional class 
satisfies the
relation
 $$
\bigl\langle\, \mu(h)^a\nu^{b+2},\ [M_{d+4}]\,\bigr\rangle
\ =\ 4\,\bigl\langle\, \mu(h)^a\nu^b,\ [M_d]\,\bigr\rangle,
\tag 1
 $$
for all $d$ and all $h\in H_2(X)$.  This condition is 
known to hold for a
large class of spaces, including the manifolds underlying 
the
simply connected elliptic surfaces, complete intersections 
and various
branched covers, as well as many ``fake'' 4-manifolds made 
from these by
surgeries and gluings.  This list probably contains all 
examples for
which any calculations have yet been made, and it is not 
impossible that (1)
is valid for all simply connected 4-manifolds.

  Whenever (1) holds, the evaluation of $\mu(h)^a\nu^b$ 
can be reduced
either to the case $b=0$, which gives the original $q_d$, 
or to the case
$b=1$. As a matter of notation, it is then convenient to 
incorporate the
case $b=1$ by introducing
invariants $q_d$ for all $d\equiv {1\over2}(b^++1)$ mod 
$2$ (rather than
just mod $4$) by defining   
 $$
2\, q_{d-2}(h)\ =\ \bigl\langle\,
\mu(h)^{d-2}\nu,\ [M_d]\,\bigr\rangle.
 $$
If $d$ is not equal to ${1\over2}(b^++1)$ mod $2$, or if 
$d$ is negative, we
set $q_d=0$.  We then combine all the polynomials $q_d$ 
into one analytic
function or formal power series $q : H_2(X,\Bbb R) \to 
\Bbb R$ by defining
 $$
q(h)\ =\ \sum_d q_d(h)/d!
 $$
Our main result is summarized in the following two theorems.

\proclaim {Theorem 1}
If $X$ is a simply connected \RM4-manifold of simple type, 
then there exist
finitely many cohomology classes $K_1$, \dots, $K_p\in
H^2(X,\Bbb Z)$ and non-zero rational numbers $a_1$, \dots, 
$a_p$ such that
 $$
q\ =\ \exp\left({Q\over2}\right)\,\sum_{s=1}^p a_s e^{K_s}
 $$
as analytic functions on $H_2(X,\Bbb R)$.  Here $Q$ is the 
intersection
form, regarded as a quadratic function.  Each of the 
\RM{``}basic classes\RM{''}
 $K_s$
is an integral lift of $w_2(X)$.
\endproclaim

\proclaim {Theorem 2}
Let $X$ be again a simply connected \RM4-manifold of 
simple type, and let
$\{K_s\}$ be the set of basic classes given by 
Theorem~\RM1.  If $\Sigma$
is any smoothly embedded, essential connected surface in 
$X$ with normal
bundle of non-negative degree, then the genus of $\Sigma$ 
satisfies the
lower bound 
 $$
2g-2 \ \ge\ \Sigma\cdot\Sigma\ +\ \max_s\, K_s\cdot \Sigma.
 $$
\endproclaim

We make some remarks about these results.  First of all, 
since $q$ is always
either an even or an odd function depending on the parity of
${1\over2}(b^++1)$, the non-zero basic classes come in 
pairs: if $K$ is a
basic class then so is $-K$, and the sum of exponentials 
which appears in
Theorem~1 is a sum of hyperbolic cosines in the even case 
or hyperbolic
sines in the odd case.  It follows too that the function
 $$
J = \max_s K_s
 $$
which appears in Theorem 2 is non-negative and even; it 
can be thought of as
defining a piecewise linear semi-norm on $H_2(X)$, and the 
inequality
 $$
2g-2 \ge Q + J
 $$
of the theorem makes it reminiscent of the Thurston norm 
on the homology of
a 3-manifold \cite{14}.  As we said in the introduction, 
the inequality is
also reminiscent of the adjunction formula for the genus 
of a smooth complex
curve: we choose the notation $K$ for the basic classes 
because we think of
them as generalizations of the canonical class of a 
complex surface.  Note
that the function $J$ can also be extracted from the 
asymptotics of $q$.  As
$h\to\infty$ in $H_2(X,\Bbb R)$, we have
 $$
\log q(h) = Q(h)/2 + J(h) + O(1).
 $$

\heading3. Examples\endheading

The first example is the $K3$ surface, whose invariants 
are known to be
given by the formula
 $$
q_{2i} = {(2i)!\over 2^ii!} Q^i
 $$
for all even degrees $2i$.  The analytic function $q$ 
therefore has the
expression 
 $$
q = \exp\left({Q\over2}\right).
 $$
Thus the only basic class for $K3$ is the zero class.  
Consider next the
complex surface $X$ formed as a double cover of $\CPtwo$, 
branched along a
smooth octic curve.  The invariants of $X$ are known to be 
polynomials in
the canonical class $K_{\!X}$ and the intersection form 
$Q$; so the basic
classes $K$ must all be multiples of $K_{\!X}$---integer 
multiples, in fact,
since $K_{\!X}$ is primitive.  Because of the last clause 
of Theorem 1, only
odd integer multiples can occur, as $w_2(X)$ is non-zero.  
Now $X$ contains
smooth complex curves $\Sigma$ whose genera are given by 
the adjunction
formula $2g-2 = \Sigma\cdot\Sigma + K_{\!X}\cdot\Sigma$, 
so according to
Theorem~2 we have
 $$
K(\Sigma) \le K_{\!X}(\Sigma),
 $$
for all basic classes $K$.  Putting these facts together, 
we see that the
only basic classes are $\pm K_{\!X}$.  Since the function 
$q$ is even in
this case and since one coefficient of one polynomial was 
calculated in
\cite{4}, we have enough information to conclude that
 $$
q = 2\,\exp\left({Q\over2}\right)\,\cosh K_{\!X}.  \tag 2
 $$
This gives us the entire invariant for the complex surface 
$X$.

Elliptic surfaces provide other examples where Theorem~1 
can be combined
with previous calculations of particular coefficients 
\cite{5} to yield a
complete answer.  For a simply connected minimal elliptic 
surface with no
multiple fibres, the invariant can be shown to be
 $$
q = \exp\left({Q\over2}\right)\,(\sinh F)^{p_g-1},  \tag 3
 $$
where $p_g$ is the geometric genus (which should be 
positive) and $F$ is the
cohomology class dual to the generic elliptic fibre.  The 
basic classes in
this case are $nF$ for integers $n$ in the range $|n| \le 
p_g-1$ satisfying
$n \equiv p_g-1$ mod $2$.  Note that $(p_g-1) F$ is the 
canonical class of
the surface.

The effect of forming a connected sum with $\CPtwobar$ in 
these examples can
also be calculated.  If $X$ is $K3\#\CPtwobar$, for 
example, we have
 $$
q = \exp\left({Q\over2}\right)\,\cosh E,  \tag 4
 $$
where $E$ is dual to the generator of $\CPtwobar$ and $Q$ 
denotes the
intersection form of $X$ (rather than the intersection 
form of $K3$).  The
basic classes are therefore $\pm E$.  The formula has the 
same shape in
examples (2) and (3) also: thus if $\tilde X = 
X\#\CPtwobar$ where $X$ is
the double cover of $\CPtwo$ branched over a smooth octic, 
then the basic
classes are $\pm K_{\!X} \pm E$ and the invariant is
 $$
q = 2\,\exp\left({Q\over2}\right)\,\cosh K_{\!X}\,\cosh E.
 $$

\heading 4. Structure of the proof\endheading

It is convenient to introduce the function $C = 
\exp(-Q/2)q$ on $H_2(X,\Bbb
R)$ and to separate it into its homogeneous parts:
 $$
C(h) = \sum_d C_d(h)/d!
 $$
Thus $C_d$ is a polynomial of degree $d$ on $H_2(X,\Bbb 
R)$, and it is
related to $q_d$ by
 $$
C_d = \sum_{i=0}^{[d/2]} {(-1)^id!\over(d-2i)!\,i!\,2^i}\, 
Q^iq_{d-2i}.  \tag 5
 $$
The content of Theorem 1 is that $C(h)$ is a linear 
combination of
exponentials.  If we focus on a particular class $S$ in 
the positive cone in
$H_2(X,\Bbb Z)$, then the central thrust of the result is 
contained in the
following proposition.

\proclaim{Proposition 3}
For any $S\in H_2(X,\Bbb Z)$ with $Q(S)$ positive, there 
is an expression
 $$
C_d(S) = \sum_{s\in\Bbb Z} \alpha_s s^d,
 $$
valid for all $d\ge0$, with only finitely many non-zero 
terms.  If $S$ is
represented by an embedded surface of genus $g$, then 
$\alpha_s$ is zero for
$|s| > 2g-2-Q(S)$.  Further, $\alpha_s$ is non-zero only 
when $s\equiv Q(S)$
$\roman{mod}\, 2$.
\endproclaim

We can view this as saying that the sequence $\{C_d(S)\}$ 
($d\in\Bbb N$)
satisfies a finite-order linear recurrence relation with 
integer roots.
These recurrence relations are the universal relations 
mentioned in the
introduction. 

Let $(X,\Sigma)$ be a pair consisting of an oriented 
4-manifold of simple
type and an oriented embedded surface. In \cite{9} and 
\cite{10} it was shown
how to associate to $(X,\Sigma)$ a family of polynomial 
invariants
$q_{k,l}$, generalizing the ordinary polynomial invariants 
of $X$.  Their
definition followed the definition of the usual 
polynomials \cite{3} but
used moduli spaces of instantons with a specified 
singularity along
$\Sigma$. The degree of $q_{k,l}$ is a function of $k$ and 
$l$ and the
topology of the pair, and like the degree of $q_d$, its 
parity is
constrained. We are interested only in the case that 
$q_{k,l}$ has degree
zero, so for any given $l$, we define $r_l$ to be 
$q_{k,l}$ if we can
find $k$ such that the degree is zero.  Otherwise, we 
define $r_l$ to
be zero.  Thus for each $l\in\Bbb Z$, we have an 
integer-valued invariant
$r_l$ for pairs $(X,\Sigma)$.

One of the main results of \cite{10} is that if $\Sigma$ 
has odd genus and
positive square, and if we set $l_{\roman o} = (g-1)/2$, 
then the invariant 
$r_{l_{\roman o}}$ can be expressed in terms of the 
ordinary polynomials: 
 $$
r_{l_{\roman o}} = 2^gq_0. \tag 6
 $$
To obtain the recurrence relation, we first generalize (6) 
to express 
$r_{l_{\roman o}-p}$ in terms of the $q_d$.  If $p<0$ then 
this invariant vanishes
\cite{10}, and for $p>0$ we establish a universal formula
 $$
r_{l_{\roman o}-p} =  A^{p,0}q_{2p}(\Sigma) + 
A^{p,2}q_{2p-2}(\Sigma)  + \cdots\ ,
\tag 7
 $$
where $A^{p,2i}$ is a quantity depending on the genus and 
self-intersection
number of $\Sigma$.  The leading term $A^{p,0}$ is 
non-zero.  This formula
is not proved in complete generality, but it is shown to 
hold when the genus
and self-intersection number of $\Sigma$ are large 
compared to $p$.

The invariants $r_l$ satisfy rather simple relations which 
are quite
easy to prove \cite{10, 8}. In particular, if the homology 
class of $\Sigma$ is
divisible by $2$, then we have
 $$
r_{l_{\roman o}-p} = \pm r_{l_{\roman o}-p'} \qquad 
\text{where} \qquad
p+p' = {\textstyle{1\over2}}(2g-2 - \Sigma\cdot\Sigma).
 $$
Using the formulae (7), we then obtain linear relations 
amongst the values
of $q_d(\Sigma)$.  When expressed in terms of the values 
of $C_d(\Sigma)$
using the formula (5), these linear relations take the 
form of a recurrence
relation with integer roots. The degree of the recurrence 
relation depends
on the genus of $\Sigma$ on account of the relationship 
between $p$ and
$p'$. The argument only establishes the validity of the 
recurrence relation
when $d$ is small compared to the genus and 
self-intersection number,
because the formula (7) is not proved for large $p$.  So 
given a homology
class $S$ of positive square, we apply the argument to 
embedded surfaces
$\Sigma$ representing various large multiples of $S$ to 
obtain eventually a
recurrence relation on $C_d(S)$ valid for all $d$.

The proof is made more complicated by the fact that we are 
unable directly
to calculate the coefficients $A^{p,2i}$ in the formula 
(7).  Instead, we
establish the nature of the resulting linear relations by 
an indirect
argument based on known examples.  The main input is our 
knowledge of the
invariants for $K3$ and some coefficients for elliptic 
surfaces without
multiple fibres \cite{12}.  An important stepping-stone is 
the proof of the
formula (4) for the invariants of $K3\#\CPtwobar$.  
Details of the proof, as
well as some more detailed applications, will appear in a 
later paper.

\heading5. Questions\endheading

It is natural to ask whether the basic classes $K_s$ can 
be shown to satisfy
any other constraints relating perhaps to the homotopy 
type of $X$.  The
canonical class of a complex surface $X$ satisfies
 $$
K_{\!X}^2 = 2\chi + 3\sigma
 $$
where $\chi$ is the Euler number and $\sigma$ is the 
signature of $X$.  In
the few examples which we know, all the basic classes 
$K_s$ satisfy this
same constraint, so it is possible that it is a general 
property.  There is
little to go on at the moment. 

One must also ask to what extent the inequality for the 
genus in Theorem~2
is sharp.  The function $J$ on $H_2(X,\Bbb Z)$ which 
appears there is
essentially the same as the function $J$ which was defined 
in section 5(iii)
of \cite{8}, at least on the positive cone. In the case 
that $X$ is a
complex surface whose invariants are polynomials in $Q$ 
and $K_{\!X}$, one
can show that $J(S) = |K_{\!X}\cdot S|$ by combining the 
results of \cite{7}
with the material of this paper.  This means showing that 
$K_{\!X}$ is one
of the basic classes.

The question of whether Theorem~2 is sharp might be 
examined in connection
with a relationship between $J$ and the Thurston norm.  
Let $Y$ be a
3-manifold with non-trivial homology $H_2(Y,\Bbb R)$, and 
consider the flat
connections over $Y$.  If we use an $\SO(3)$ bundle with 
suitably chosen
Stiefel-Whitney class, we may arrange that there are no 
reducible
connections; and under these circumstances it is possible 
to define a Floer
homology group $\text{\it HF\/}(Y)$.  The construction of 
the polynomial
invariants gives, in this situation, a linear map $\Phi$ 
from $H_2(Y,\Bbb
R)$ to the ring of endomorphisms of $\text{\it HF\/}(Y)$.  
All elements of
the image of $\Phi$ commute \cite{2}, so they have 
simultaneous eigenvalues
which can be regarded as defining $r$ linear functions on 
$H_2(Y,\Bbb R)$,
where $r$ is the rank of the Floer homology: these are 
close cousins of the
basic classes in this paper.  On the basis of the results 
we have described,
one might speculate that the linear functions defined by 
the eigenvalues are
integral and that their supremum defines a semi-norm $J_Y$ 
on $H_2(Y,\Bbb
R)$ which provides a lower bound for the Thurston norm 
$x$.  It would then
be interesting to know if there is any significant class 
of 3-manifolds for
which $x$ and $J_Y$ can be shown to be equal.

The role of the simple type condition should also be 
clarified.  If there
are simply connected 4-manifolds which are not of simple 
type, it is
possible that the statement of Theorem~1 would need only 
minor modification
to cover the general case.  Perhaps the constant 
coefficients $a_s$ which
appear there should be replaced by polynomial functions.

We close with a conjecture about the invariants of 
elliptic surfaces, rather
less speculative than the remarks above.  Let
$X$ be a regular elliptic surface with geometric genus 
$p_g\ge1$
and $r$ multiple fibres of multiplicities $m_1$, \dots, 
$m_r$.  We want $X$
to be simply connected, which means that we must have 
$r\le2$ and that in the
case $r=2$ the two multiplicities must be coprime. Then we
conjecture that
 $$
q \ =\ \exp\left({Q\over2}\right)\,
{(\sinh F)^{p_g-1+r}\over \prod_i \sinh(F/m_i)}\,,
 $$
where $F$ is the class dual to the generic fibre. For 
$r=0$ this is the
formula given in section~3. For $r=1$ the formula is 
correct in a few
small cases, such as when $p_g=1$ and $m_1\le4$, or 
$p_g=2$ or $3$ and
$m_1=2$.  For $r=2$ and all $p_g$, the formula is in 
agreement with the
known expressions for the first two non-zero coefficients 
of $q_d$ in its
expansion as a polynomial in $Q$ and $F$ \cite{12}.  The 
results of this
paper will probably extend without change to the 
non-simply connected case,
as long as $H_1(X,\Bbb R)=0$, so the restriction $r\le2$ 
should not be
regarded as essential.


\Refs

\ref\no 1
\by S. Bauer
\paper	Diffeomorphism types of elliptic surfaces with 
$p_g=1$
\paperinfo Warwick preprint, 1992
\endref

\ref\no 2
\by	P. J. Braam and S. K. Donaldson
\paper	Fukaya-Floer homology and gluing formulae for 
polynomial invariants
\paperinfo preprint
\endref

\ref\no 3
\by	S. K. Donaldson
\paper	Polynomial invariants for smooth four-manifolds
\jour	Topology
\vol	29 \yr 1990 \pages 257--315
\endref

\ref\no 4
\by	S. K. Donaldson and P. B. Kronheimer
\book	The geometry of four-manifolds
\publ	Oxford Univ. Press \publaddr Oxford \yr 1990
\endref

\ref\no 5
\by	R. Friedman and J. W. Morgan 
\paper	Complex versus differentiable classification of 
algebraic
surfaces
\jour	Topology Appl.
\vol	32 \yr 1989 \pages 135--139
\endref

\ref\no 6
\bysame
\book	Smooth four-manifolds and complex surfaces 
\publ	Springer-Verlag \publaddr New York\toappear
\endref

\ref\no 7
\by	P. B. Kronheimer
\paper	The genus-minimizing property of algebraic curves
\jour	Bull. Amer. Math. Soc. (N.S.) \vol 29 \yr 1993 
\pages 63--9
\endref

\ref\no 8
\bysame
\paper	An obstruction to removing intersection points in 
immersed surfaces 
\jour Topology \toappear
\endref

\ref\no 9
\by	P. B. Kronheimer and T. S. Mrowka
\paper	Gauge theory for embedded surfaces \rom{I}
\jour	Topology 
\vol 32 \yr1993 \pages 773--826
\endref

\ref\no 10
\bysame
\paper	Gauge theory for embedded surfaces \rom{II}
\jour	Topology \toappear
\endref

\ref\no 11
\by	J. W. Morgan and T. S. Mrowka 
\paper  A note on Donaldson\RM's polynomial invariants
\jour   Internat. Math. Res. Notices \vol 10 \yr 1992 
\pages 223--230
\endref

\ref\no 12
\bysame	
\paper	On the diffeomorphism classification of regular 
elliptic surfaces 
\jour	Internat. Math. Res. Notices \toappear
\endref

\ref\no 13
\by	J. W. Morgan and K. G. O'Grady
\paper	The smooth classification of fake K\RM{3'}s and 
similar surfaces
\paperinfo preprint, 1992
\endref

\ref\no 14
\by	W. P. Thurston
\paper	A norm for the homology of \RM3-manifolds
\inbook	Mem. Amer. Math. Soc., vol. 59 
\publ Amer. Math. Soc. \publaddr Providence, RI
\yr 1986  \pages 99--130
\endref
\endRefs
\enddocument